# NAVIGATING TO THE MOON ALONG LOW-ENERGY TRANSFERS

M. Vetrisano[1], W. Van der Weg[2], M. Vasile[3]

Department of Mechanical & Aerospace Engineering, University of Strathclyde
75 Montrose Street G1 1XJ, Glasgow, UK, t:0141 548 2083, f:0141 552 5105, e-mails
{massimo.vetrisano, willem.van-der-weg, massimiliano.vasile}@strath.ac.uk

**Abstract**. This paper presents a navigation strategy to fly to the Moon along a Weak Stability Boundary transfer trajectory. A particular strategy is devised to ensure capture into an uncontrolled relatively stable orbit at the Moon. Both uncertainty in the orbit determination process and in the control of the thrust vector are included in the navigation analysis. The orbit determination process is based on the definition of an optimal filtering technique that is able to meet accuracy requirements at an acceptable computational cost. Three sequential filtering techniques are analysed: an extended Kalman filter, an Unscented Kalman filter and a Kalman filter based on high order expansions. The analysis shows that only the unscented Kalman filter meets the accuracy requirements at an acceptable computational cost. This paper demonstrates lunar weak capture for all trajectories within a capture corridor defined by all the trajectories in the neighbourhood of the nominal one, in state space. A minimum $\Delta v$ strategy is presented to extend the lifetime of the spacecraft around the Moon. The orbit determination and navigation strategies are applied to the case of the European Student Moon Orbiter.

## 1 INTRODUCTION

Low cost missions to the Moon generally impose limits on the control over launch opportunities (e.g. piggy-backing on a scheduled commercial launch) and on the qualification of a number of components, such as the attitude control and propulsion systems, which translates into an increased uncertainty of the overall mission. Low-energy transfer solutions, such as a Weak Stability Boundary (WSB) transfer, can be used to increase transfer flexibility and reduce the transfer Δv demand. However, the chaotic nature of low-energy transfers and the low confidence level in key spacecraft components such as propulsion and attitude control impose the use of robust navigation strategies based on a reliable estimation of the trajectory. For these types of missions, orbit uncertainty quantification and propagation plays a major role in orbit determination and correction manoeuvre design.

Furthermore, such missions require the definition of particular navigation strategies that are able to guarantee capture at the Moon under contingency situations at minimal Δv cost. The objective of classical trajectory navigation is to follow the reference trajectory while minimizing the fuel consumption. In this paper, a different approach is proposed. The idea is to maintain the spacecraft within a so called capture corridor to achieve satisfactory orbit insertion at the Moon into a lunar orbit suited to payload's requirements. This approach allows for the following of a trajectory which is considerably far from the reference one. Because the reference and the actual trajectory are expected to be very different, it is desirable to use a filtering technique that is able to deal with nonlinear system dynamics and combines precision and process velocity. In this manner data processing can be performed either in real time while measurements are being collected or at a later time when the full set of measurements is available. In this paper, three filters derived from the Kalman filter were put to the test: the extended Kalman filter (Battin, 1970), the unscented Kalman filter (Julier et al., 1995), and the higher order semi-analytic extended Kalman filters developed by Park and Scheeres, 2006a-b. These types of filter are all suitable for nonlinear systems but in the framework of the case study presented in this work, it will be shown that the unscented filter provides an advantage compared to the other two.

This paper investigates the application of this navigation and orbit determination methodology to the case of the WSB transfer to the Moon planned for the European Student Moon Orbiter. In this investigation, the possibility for the main engine to fail before the lunar orbit insertion is taken into account in order to evaluate the robustness of the whole navigation strategy. A contingency analysis is performed on trajectories where the characteristic energy $C_3$ at lunar arrival is negative, as the spacecraft can be captured by the Moon for a relatively long period. The attitude thrusters are considered able to produce small correction manoeuvres and extend the orbit lifetime around the Moon. The results show that it is possible to extend the spacecraft lifetime around the Moon by applying small Δv manoeuvres. By keeping the spacecraft within the corridor it is possible to recover the spacecraft and increase the mission probability of success.

---

[1] PhD Candidate, Space Advanced Research Team, Advanced Space Concepts Laboratory
[2] PhD Candidate, Space Advanced Research Team, Advanced Space Concepts Laboratory
[3] Reader, Space Advanced Research Team, Advanced Space Concepts Laboratory

## 2 CASE STUDY

The European Student Moon Orbiter (ESMO) is the fourth small satellite mission within ESA's Education Satellite Program, and is currently the only scheduled ESA mission to the Moon for the next four years (Walker and Cross, 2010). The spacecraft is designed, built, and operated entirely by students, under the guidance of Surrey Satellite Technology Limited. ESMO is currently scheduled for launch in 2014 – 2015 as a secondary payload that will be inserted into a Geostationary Transfer Orbit (GTO). The shared launch allows for little to no control over exact launch date and time, and so a substantial amount of analysis must be performed in order to cover the two years of possible launch opportunities. The scientific payload consists of a narrow angle camera to take high resolution images of the lunar surface, where particular attention is paid to the observation of the South Pole.

The spacecraft will use its chemical propulsion system (Croisard et al., 2009) to transfer itself from GTO to its lunar operational orbit at the Moon using a WSB transfer (Koon et al., 2001), as shown in Figure 1. This type of transfer has been selected due to its associated propellant saving, and to cope consistently with a variety of injection conditions resulting from the shared launch. Limitations of the propulsion system imply that the departure conditions cannot be met in a single manoeuvre. Departure from the Earth (i.e. injection of the spacecraft into its lunar bound WSB transfer) is performed by employing a number of apogee raising manoeuvres before using a final manoeuvre to begin the transfer to the Moon.

The critical design driver for ESMO is the mission cost. In addition to the aforementioned shared launch and propellant saving WSB transfer, further cost saving measures are also implemented. Commercial off the shelf parts and flight spares are extensively used in order to reduce cost. The use of a highly eccentric frozen polar orbit (Gibbings et al., 2010) serves to reduce the propellant necessary to inject into the final lunar orbit, and to reduce the propellant necessary to maintain the operational orbit. Finally, a single ground station reduces the operating cost of the ground segment. The cost saving comes with an increased uncertainty on a number of critical aspects of the transfer.

### 2.1 *ESMO Earth-Moon Transfer*

ESMO will use a type of low-energy transfer first proposed by Edward Belbruno, in 1987 (Belbruno, 1987). This type of transfer employs the combined gravity force of the Sun, Earth, and Moon advantageously in order to change the orbit plane, and to raise the perigee of the orbit from the Earth up to the Moon. In an inertial geocentric reference frame, the motion of the spacecraft is governed by the following set of differential equations (Bate et al., 1971):

$$\dot{\mathbf{r}} = \mathbf{v}$$
$$\dot{\mathbf{v}} = -\frac{\mu_E}{r^3}\mathbf{r} - \mu_S\left(\frac{\mathbf{r}_{Ssc}}{r_{Ssc}^3} - \frac{\mathbf{r}_{SE}}{r_{SE}^3}\right) - \mu_M\left(\frac{\mathbf{r}_{Msc}}{r_{Msc}^3} - \frac{\mathbf{r}_{ME}}{r_{ME}^3}\right) \quad [1]$$

where **r** and **v** are respectively the position and velocity vectors of the spacecraft with respect to the Earth-in the J2000 inertial reference frame, $\mathbf{r}_{Ssc}$ and $\mathbf{r}_{SE}$ are the Sun-spacecraft and Sun-Earth vectors, $\mathbf{r}_{Ssc}$ and $\mathbf{r}_{SE}$ are the Moon-spacecraft and Moon-Earth vectors, $\mu_E$, $\mu_S$ and $\mu_M$ are the planetary constants of Earth, Sun and Moon respectively. The position of Sun and Moon with respect to the Earth and the spacecraft are calculated using analytical ephemeris (Vallado 2004), accounting for secular variations in the orbital elements of both the Earth and the Moon. An analytical model was used to describe the secular variation of the angles between the Earth-equatorial and the Moon-equatorial reference frame. Eqs. [1] describe the perturbed motion of the spacecraft where the second and third terms in the variation of the velocity vector are the gravity perturbation of the Sun and the Moon respectively.

If compared to a more traditional Hohmann transfer, this kind of trajectory increases the transfer time by approximately 3 months. Figure 1 shows a typical WSB transfer for ESMO. In principle, 3 manoeuvres would suffice to perform a complete transfer − one manoeuvre at the Earth, followed by a single course correction, and finally a manoeuvre to inject into orbit at the Moon. The transfer in Figure 1 was generated by forward propagation of Eqs. [1] from time $t_0$ to time $t_{WSB}$ with initial conditions from a standard Ariane 5 GTO and by backwards propagation of Eqs. [1] from time $t_f$ up to $t_{WSB}$. The course correction manoeuvre is performed at time $t_{WSB}$ in correspondence of the WSB point, which represents the connecting point between the forward propagated arc from the Earth and the backward propagated arc from the Moon (c.f. Figure 1 where the black and red arcs meet). The location of the point is not predefined but is an outcome of the transfer design optimisation.

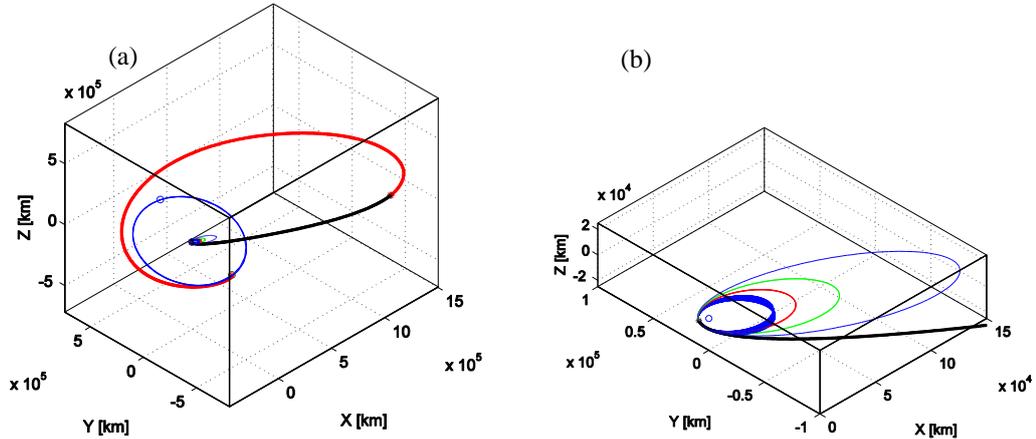

Figure 1:Typical trajectory for WSB transfer in the Earth Centred Inertial frame. The blue line represents the Moon's orbit around the Earth, the black-red line is the WSB transfer from Earth to Moon (a), which is preceded by the apogee raising strategy (b) close to the Earth.

The total cost in $\Delta v$ to achieve translunar injection is in the region of 750 m/s. A small manoeuvre, with a magnitude that can range between 0 and 100m/s, is performed at time $t_{WSB}$, allowing the spacecraft to return back toward the Earth-Moon system, but now such that the spacecraft can be inserted into lunar orbit using a final manoeuvre of roughly 100 m/s.

The launch may occur at any date within the launch window (2014-2015) therefore a full database of possible transfer was populated with over 300,000 transfer trajectories (Van der Weg and Vasile, 2011). A useful parameter to characterize each transfer trajectory at the Moon is represented the $C_3$, which is defined as $-\mu_M / a$ (where $a$ is the semi-major axis of the incoming hyperbolic trajectory at the Moon). This is derived from the definition of the two-body energy for a massless particle orbiting the central body. The orbit is open when the $C_3$ value is positive, and closed when it is negative. Figure 2 shows the $C_3$ values for the year 2014.

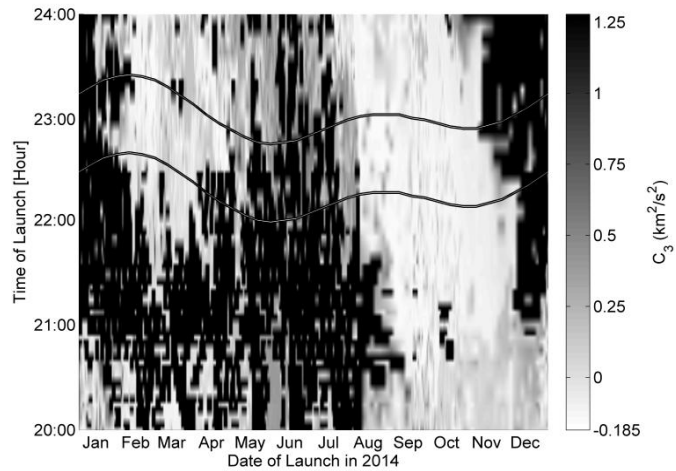

Figure 2: Plot showing $C_3$ value at lunar arrival as a function of departure date and launch hour for the Ariane 5 launch vehicle. The 'wavy' lines in white superposed on top of the plot indicate the standard daily launch window for a commercial Ariane 5 GTO launch.

For every date and launch time between 20:00 and 24:00 UTC a valid transfer solution is sought that satisfies mission constraints (e.g. the desired lunar orbit parameters). Every valid transfer possesses a value of $C_3$ at lunar arrival. This is shown in Figure 2, which shows the $C_3$ values for the launch conditions in 2014 using an Ariane 5 launch vehicle. Areas in black indicate combinations of date and time that did not result in the generation of a valid transfer.

The value of the $C_3$ at lunar arrival varies on the location of the WSB point and the manoeuvre magnitude. However, due to the shifting Earth-Sun-Moon geometry over time and the selected set of orbital parameters (for both targeted lunar orbit and departure GTO) the transfer changes. The location of the WSB point and the manoeuvre magnitudes change, thus the $C_3$ at lunar arrival will change (this is a by-product of the optimization search for a propellant optimal transfer). The plot not only shows the influence of a different departure date (Earth-Sun-Moon geometry shifts) throughout the year, but that there is an additional influence from the launch time for a particular date. The time of launch for each day influences the orientation of the GTO with respect to the Sun (i.e.: a change in the time of launch changes the angle between the Earth-Sun vector and line of apsides of the resulting GTO). This orientation also has an influence on the location of the manoeuvre point and magnitude as the optimization will select a different WSB point and manoeuvre magnitude.

Many trajectories in Figure 2 present a negative $C_3$ value at Lunar Orbit Insertion (LOI). In general a negative $C_3$ is a favourable condition as the spacecraft can be temporarily captured by the Moon without requiring manoeuvres. In this way the spacecraft can experience a so called weak capture (Makó et al., 2010), which leads the spacecraft to temporarily orbit around the Moon without having performed the LOI manoeuvre. In a two-body system the energy would remain constant, and an external force would be necessary to change the energy state. However, the gravity of third bodies such as the Earth and the Sun (as well as other perturbations such as those from the lunar gravity field) has a substantial impact on the duration of time for which the spacecraft is weakly captured. Therefore, a negative $C_3$ is not sufficient to guarantee that the spacecraft is captured for a satisfactory period of time. Starting from the solutions with negative $C_3$, Section 6 will investigate how to improve the stability and lifetime of the weakly captured orbits when the main engine is assumed to fail before lunar insertion.

## 3 ORBIT DETERMINATION

The Orbit Determination, OD, is performed by processing raw measurements coming from the only available ground station of Raisting, located in Germany, whose coordinates are listed in Table 1.

Table 1: Raisting ground station location.

| Latitude λ [deg] | Longitude υ [deg] | Altitude above mean sea level [m] | Minimum Elevation angle [deg] |
|---|---|---|---|
| 47.90221 | 11.11579 | 553 | 1 |

The set of measurements includes range and range rate $\rho$ and $\dot{\rho}$, from the ground station, plus the pointing angles $A$, $E$ (respectively azimuth and elevation). Since the actual position of the spacecraft is given in the Earth Centred Inertial (ECI) reference frame, it is necessary to write the state of ESMO as it was seen in the local South East Zenith (SEZ) reference frame of the ground station as shown in Figure 3.

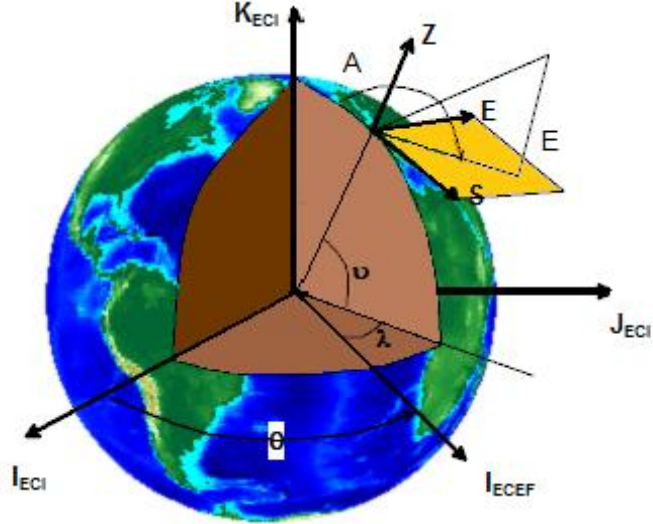

Figure 3: Reference frames (Vallado, 2000)

The range in the ECI reference frame is given by the difference in position of the spacecraft and ground station location (Montenbruck, 2000):

$$\boldsymbol{\rho}_{ECI} = \mathbf{r}_{ECI} - \mathbf{r}_{site-ECI} \qquad [2]$$

where $\mathbf{r}_{site-ECI}$ is the ECI position of the ground station. The range and velocity vectors in the SEZ frame are given by the following transformations:

$$\begin{aligned}\boldsymbol{\rho}_{SEZ} &= \mathbf{A}_{SEZ-ECEF}\,\mathbf{A}^{T}_{ECI-ECEF}\,\boldsymbol{\rho}_{ECI} \\ \dot{\boldsymbol{\rho}}_{SEZ} &= \mathbf{A}_{SEZ-ECEF}\,\mathbf{A}^{T}_{ECI-ECEF}\,\mathbf{v}_{ECI}\end{aligned} \qquad [3]$$

where $\mathbf{A}_{SEZ\text{-}ECEF}$ and $\mathbf{A}_{ECI\text{-}ECEF}$ are respectively the transformation matrix from the Earth Centred Earth fixed (ECEF) reference frame to SEZ and the transformation matrix from ECEF to ECI:

$$\mathbf{A}_{ECI-ECEF} = \begin{pmatrix} \cos\theta & -\sin\theta & 0 \\ \sin\theta & \cos\theta & 0 \\ 0 & 0 & 1 \end{pmatrix}$$

$$\mathbf{A}_{SEZ-ECEF} = \begin{pmatrix} \sin\upsilon\cos\lambda & -\sin\upsilon\sin\lambda & -\cos\lambda \\ -\sin\lambda & \cos\lambda & 0 \\ \cos\upsilon\cos\lambda & \cos\upsilon\sin\lambda & \sin\upsilon \end{pmatrix} \qquad [4]$$

λ and υ are the latitude and longitude of the ground station; $\theta = 280.4606 + \Omega t$ is the rotation angle between the ECI and ECEF reference frame about the z-axis; Ω is the Earth's angular speed expressed in deg/day and $t$ is the time expressed in MJD2000. Finally the set of simulated measurements is obtained from the SEZ position and velocity:

$$\begin{aligned}\rho &= |\boldsymbol{\rho}_{SEZ}| \\ E &= \sin^{-1}\frac{\rho_Z}{\rho} \\ A &= \sin^{-1}\frac{\rho_E}{\sqrt{\rho_S^2 + \rho_E^2}} \\ \dot{\rho} &= \frac{\dot{\boldsymbol{\rho}}_{SEZ}\cdot\boldsymbol{\rho}_{SEZ}}{\rho}\end{aligned} \qquad [5]$$

being $\rho_S$, $\rho_E$, $\rho_Z$ the components of the range in the SEZ frame.

The actual measurements were simulated by perturbing the nominal ones defined in Eq.[5] with a random noise with normal distribution. Ionospheric and tropospheric refractions were not taken into account since their effects are expected to be corrected at ground station level before the information is provided to flight dynamics for orbit determination. The orbit state vector is obtained by means of a suitable filtering technique able to deal with nonlinearities in the system dynamics.

## 3.1 *Filtering Techniques*

The dynamics and measurement model used in the filtering takes the following form (Montenbruck and Gill, 2000):

$$\dot{\mathbf{x}}(t) = \begin{bmatrix} \dot{\mathbf{r}} \\ \dot{\mathbf{v}} \end{bmatrix} = \begin{bmatrix} \mathbf{v} \\ \mathbf{f}(\mathbf{x}(t),t) \end{bmatrix} + \begin{bmatrix} \mathbf{0} \\ \mathbf{a}_{Solar} + \mathbf{a}_w \end{bmatrix} \quad [6]$$

$$\mathbf{y}_k = \mathbf{h}(\mathbf{x}_k) + \mathbf{v}_k$$

where $\mathbf{f}(\mathbf{x}(t),t)$ is the set of nonlinear continuous-time equations in the second of Eqs. [1], $\mathbf{a}_w$ is the random noise acceleration produced by the attitude system, $\mathbf{a}_{Solar}$ is the solar pressure modelled as follows:

$$\mathbf{a}_{Solar} = -C_R F_{Solar}(r_{Ssc}) \frac{A}{M(t)} \frac{\mathbf{r}_{Ssc}}{r_{Ssc}} \quad [7]$$

where $A$, $M$ are the satellite cross section area and mass respectively, $C_R$ the reflectivity coefficient and $F_{Solar}$ is the solar flux. In Eqs. [6] the measurement model $h(x_k)$ corresponds to the set of four equations in Eqs.[5], giving the values of range, elevation, azimuth and range-rate, and $\mathbf{y}_k$ is the set of actual measurements corresponding to the state vector $\mathbf{x}_k$ at the measurement time $t_k$. The term $\mathbf{v_k}$ is the measurements noise defined as $\mathbf{v}_k = \mathbf{R}\boldsymbol{\sigma}_{normal}$, where $\boldsymbol{\sigma}_{normal} \in N(0,1)$ is a vector of random numbers taken from a normal distribution. R is the observation error covariance matrix whose components along the diagonal are the squared values of the ground station errors associated with the measurements. Since it is assumed that all the measurements are uncorrelated and independent, extra diagonal elements in R are null and the matrix is diagonal.

In the following three sequential filtering techniques are analysed to establish which approach is the most suitable to navigate a spacecraft to the Moon along a WSB transfer. The three filters are: an Extended Kalman Filter (EKF), an Unscented Kalman Filter (UKF) and a Kalman filter based on high order Taylor expansions.

### 3.1.1 Extended Kalman filter

The EKF is a well-known dynamic optimal filter which was first employed in the Apollo program (Battin and Levine, 1970). The EKF linearizes the filtering problem about the estimated state and covariance matrix. It is composed of two conceptually distinct phases: the time update and the measurements update. The time update phase consists of the propagation of the latest estimate $\hat{\mathbf{x}}_k^+$ to obtain an a-priori estimate at current epoch $\hat{\mathbf{x}}_{k+1}^-$ with the corresponding covariance matrix $\mathbf{P}_{k+1}^-$ (Montenbruck and Gill, 2000):

$$\hat{\mathbf{x}}_{k+1}^- = \mathbf{x}(t_{k+1}; \mathbf{x}(t_k) = \hat{\mathbf{x}}_k^+)$$
$$\mathbf{P}_{k+1}^- = \boldsymbol{\Phi}_{k+1} \mathbf{P}_k^+ \boldsymbol{\Phi}_{k+1}^T \quad [8]$$

where the predicted state $\hat{\mathbf{x}}_{k+1}^-$ is obtained by integrating forward Eqs. [6] starting from the latest estimate state $\hat{\mathbf{x}}_k^+$; $\boldsymbol{\Phi}_{k+1}$ is the state transition matrix (STM) coming from the linearization of the dynamic equations about the updated state:

$$\boldsymbol{\Phi}_{k+1} = \boldsymbol{\Phi}_{k+1}(t_{k+1}, t_k) = \boldsymbol{\Phi}_{k+1}(t_k, t_k) + \int_{t_k}^{t_{k+1}} \mathbf{F}(\mathbf{x},t) \boldsymbol{\Phi}_{k+1}(t,t_k) dt \approx \mathbf{I} + \mathbf{F}_k \Delta t$$

$$\mathbf{F}_k = \frac{\partial \mathbf{f}(\mathbf{x}_k^+)}{\partial \mathbf{x}_k^+} \quad [9]$$

being $\boldsymbol{\Phi}_{k+1}(t_k, t_k)$ equal to the identity matrix. The measurement update phase consists in the computation of the Kalman gain $\mathbf{K}_{k+1}$ and the state estimate $\hat{\mathbf{x}}_{k+1}^+$ and covariance matrix, $\mathbf{P}_{k+1}^+$ updates:

$$\mathbf{K}_{k+1} = \mathbf{P}_{k+1}^{-}\mathbf{H}_{k+1}^{T}\left[\mathbf{H}_{k+1}\mathbf{P}_{k}^{-}\mathbf{H}_{k}^{T} + \mathbf{R}\right]^{-1}$$
$$\hat{\mathbf{x}}_{k+1}^{+} = \hat{\mathbf{x}}_{k+1}^{-} + \mathbf{K}_{k+1}(\mathbf{y}_{k+1} - \mathbf{h}(\hat{\mathbf{x}}_{k+1}^{-}))$$
$$\mathbf{H}_{k+1} = \frac{\partial \mathbf{h}(\hat{\mathbf{x}}_{k+1}^{-})}{\partial \hat{\mathbf{x}}_{k+1}^{-}}$$
$$\mathbf{P}_{k+1}^{+} = (\mathbf{I} - \mathbf{K}_{k+1}\mathbf{H}_{k+1})\mathbf{P}_{k+1}^{-}$$
[10]

where $\mathbf{H}_{k+1}$ is the Jacobian matrix of the measurement function. The Kalman gain as reported in Eqs. [10] minimizes the a posteriori error covariance matrix. It represents a function of the relative certainty of the measurements and current state. As the gain increases the measurements are trusted more and the estimates rely less on the prediction model. On the contrary as the gain decreases, measurements tend to be ignored and the estimate rely more heavily on the prediction model.

### 3.1.2 Unscented Kalman Filter

The unscented Kalman filter, first introduced by Julier et al. in 1995, works on the premises that by using a limited set of sample, optimally chosen, it should be easier to approximate a Gaussian distribution than to approximate a nonlinear function. The UKF was shown to be preferable to the EKF in the case of nonlinear systems as the expected error in terms of mean and covariance matrix is lower, and it does not require the derivation of the Jacobian matrix (Crassidis et al., 2000). The Kalman filter updates are represented as follows:

$$\hat{\mathbf{x}}_{k}^{+} = \hat{\mathbf{x}}_{k}^{-} + \mathbf{K}_{k}\mathbf{\upsilon}_{k}$$
$$\mathbf{P}_{k}^{+} = \mathbf{P}_{k}^{-} - \mathbf{K}_{k}\mathbf{P}_{k}^{\upsilon\upsilon}\mathbf{K}_{k}^{T}$$
[11]

where $\mathbf{\upsilon}_k$ and the Kalman gain $\mathbf{K}_k$ are:

$$\mathbf{\upsilon}_{k} \equiv \mathbf{y}_{k} - \hat{\mathbf{y}}_{k}^{-} = \tilde{\mathbf{y}}_{k} - \mathbf{h}(\hat{\mathbf{x}}_{k}^{-}, k)$$
$$\mathbf{K}_{k} = \mathbf{P}_{k}^{xy}(\mathbf{P}_{k}^{\upsilon\upsilon})^{-1}$$
[12]

$\mathbf{\upsilon}_k$ is called *innovations process* and represents the difference between the current measurements and the predicted ones, $\mathbf{P}_{k}^{\upsilon\upsilon}$ is the covariance matrix of the innovations process at the sample time $t_k$. The matrix $\mathbf{P}_{k}^{xy}$ is the cross-correlation between $\hat{\mathbf{x}}_{k}^{-}$ and $\hat{\mathbf{y}}_{k}^{-}$. The approach used in the filter design requires augmenting the covariance matrix with:

$$\mathbf{P}_{k}^{a} = \begin{bmatrix} \mathbf{P}_{k}^{+} & \mathbf{P}_{k}^{x\omega} & \mathbf{P}_{k}^{xv} \\ (\mathbf{P}_{k}^{x\omega})^{T} & \mathbf{Q}_{k} & \mathbf{P}_{k}^{\omega v} \\ (\mathbf{P}_{k}^{xv})^{T} & (\mathbf{P}_{k}^{\omega v})^{T} & \mathbf{R}_{k} \end{bmatrix}$$
[13]

where $\mathbf{P}_{k}^{xw}$ is the correlation between the state error and process noise, $\mathbf{P}_{k}^{xv}$ is the correlation between the state error and measurement noise, and $\mathbf{P}_{k}^{wv}$ is the correlation between the process noise and measurement noise, which is zero in this case. If $L$ is the number of elements per column of the augmented covariance matrix $\mathbf{P}_{k}^{a}$, then a set of $2L$ samples $\mathbf{\chi}_{k}^{a}(i)$, with $i=1,..,2L$, called sigma points, are generated such that:

$$\mathbf{\chi}_{k}^{a}(i) = \mathbf{\sigma}_{k}(i) + \hat{\mathbf{x}}_{k}^{a}$$
$$\mathbf{\sigma}_{k} = \pm\eta\sqrt{\mathbf{P}_{k}^{a}}$$
$$\mathbf{\chi}_{k}^{a}(0) = \hat{\mathbf{x}}_{k}^{a}$$
[14]

where $\eta$ is a suitable scaling factor and $\hat{\mathbf{x}}_{k}^{a}$ is the augmented state :

$$\hat{\mathbf{x}}_k^a = \begin{bmatrix} \hat{\mathbf{x}}_k \\ 0_{q \times 1} \\ 0_{z \times 1} \end{bmatrix} \qquad [15]$$

$q$ is the dimension of *w*, and $z$ is the dimension of $y_k$. The sampled sigma points are then:

$$\chi_k^a(i) = \begin{bmatrix} \chi_k^x(i) \\ \chi_k^w(i) \\ \chi_k^v(i) \end{bmatrix} \qquad [16]$$

$\chi_k^x$ is the vector of the first *n* (size of $x_k$) elements of $\chi_k^a$, $\chi_k^w$ is a vector of the next *q* elements of $\chi_k^a$ and $\chi_k^v$ is the vector of the last *z* components of $\chi_k^a$. The sigma points are transformed or propagated through the system dynamics equation by the so called *unscented transformation* (UT):

$$\begin{aligned} \chi_{k+1}(i) &= \mathbf{f}(\chi_k^x(i), \chi_k^\omega(i), \mathbf{u}_k, t_k) \\ \mathbf{\mu}_{k+1}(i) &= \mathbf{h}(\chi_{k+1}^x(i), \mathbf{u}_{k+1}, \chi_{k+1}^v(i), t_{k+1}) \end{aligned} \qquad [17]$$

The predicted mean of the state vector, the covariance matrix and the mean observation can be approximated using the weighted mean and covariance of the transformed vectors:

$$\begin{aligned} \hat{\mathbf{x}}_{k+1}^- &= \sum_{i=0}^{2L} W_i^m \chi_{k+1}^x(i) \\ \mathbf{P}_{k+1}^- &= \sum_{i=0}^{2L} W_i^c [\chi_{k+1}^x(i) - \hat{\mathbf{x}}_{k+1}^-][\chi_{k+1}^x(i) - \hat{\mathbf{x}}_{k+1}^-]^T \\ \hat{\mathbf{y}}_{k+1}^- &= \sum_{i=0}^{2L} W_i^m \mathbf{\mu}_{k+1}(i) \end{aligned} \qquad [18]$$

where $W_i^m$, $W_i^c$ are suitable weighting factors. Finally the updated covariance and the cross correlation matrix are:

$$\begin{aligned} \mathbf{P}_{k+1}^{yy} &= \mathbf{P}_{k+1}^{vv} = \sum_{i=0}^{2L} W_i^{cov}[\gamma_{k+1}(i) - \hat{\mathbf{y}}_{k+1}^-][\gamma_{k+1}(i) - \hat{\mathbf{y}}_{k+1}^-]^T \\ \mathbf{P}_{k+1}^{xy} &= \sum_{i=0}^{2L} W_i^{cov}[\chi_{k+1}^x(i) - \hat{\mathbf{x}}_{k+1}^-][\gamma_{k+1}(i) - \hat{\mathbf{y}}_{k+1}^-]^T \end{aligned} \qquad [19]$$

In this way it is possible to update the filter to the next observation and prediction at time $t_{k+1}$, using the equations reported in Eq.[11].

### 3.1.3 High Order Semi-Analytic Extended Kalman Filter

In 2006, Park and Sheeres derived analytical expressions of a nonlinear trajectory solution using higher order Taylor series approach and applied the results to spacecraft application. In particular they presented a semi-analytic filtering method by implementing the State Transition Tensors (STTs) to sequentially update the state vector with contributions from each measurements. They called this nonlinear filter High-order semi-Analytic Extended Kalman filter (HAEKF), since the implementation follows the same steps of the conventional Kalman Filter. The main characteristic and advantage of using the higher order semi-analytic extended Kalman filter is that the STTs can be calculated offline prior to their usage in the filter itself. The STTs map analytically the local nonlinear motion of the spacecraft at the current epoch to the initial deviated conditions from the nominal trajectory. This section describes the fundamental aspects which this method is based on. The local spacecraft dynamics can be described by applying a Taylor series expansion about the reference nominal trajectory $\mathbf{x}_0$ for some initial deviation $\delta\mathbf{x}_0$:

$$\delta\mathbf{x}(t) = \phi(t, \mathbf{x}_0 + \delta\mathbf{x}_0; t_0) - \phi(t, \mathbf{x}_0; t_0) \qquad [20]$$

$\phi$ is the solution flow which maps the initial state at $t_0$ to $t$. The *s*-th order solution can be expressed using the Einstein summation convention:

$$\delta x^i(t) = \sum_p^s \frac{1}{p!} \phi_{(t,t_0)}^{i,\gamma_1...\gamma_p} \delta x_0^{\gamma_1} \cdots \delta x_0^{\gamma_p} \quad [21]$$

where $\gamma_1..\gamma_p \in \{1,....,n\}$ denotes the $\gamma_i$ component of the state vector corresponding to the s-th derivative, $n$ is the number of components of the state vector and:

$$\phi_{(t,t_0)}^{i,\gamma_1...\gamma_p}(t;\mathbf{x}_0;t_0) = \left.\frac{\partial^p \phi_{(t,t_0)}^i(t;\xi_0;t_0)}{\partial \xi_0^{\gamma_1} \cdots \partial \xi_0^{\gamma_p}}\right|_{\xi_0^{\gamma_j}=x_0^{\gamma_j}} \quad [22]$$

The higher-order partials of the solution define the global state transition tensors, which map the initial deviations at time $t_0$ to the deviation at time $t$. The higher order effects are included in the STTs. Note that for $s=1$, the STTs reduces to the simple state transition matrix. The differential equations up to the third order are given by the following equations (Park and Sheeres, 2006a):

$$\dot{\phi}^{i,a} = f^{i,\alpha} \phi^{\alpha,a} \quad [23]$$

$$\dot{\phi}^{i,ab} = f^{i,\alpha} \phi^{\alpha,ab} + f^{i,\alpha\beta} \phi^{\alpha,a} \phi^{\beta,b} \quad [24]$$

$$\dot{\phi}^{i,abc} = f^{i,\alpha} \phi^{\alpha,abc} + f^{i,\alpha\beta} \left( \phi^{\alpha,a} \phi^{\beta,bc} + \phi^{\alpha,ab} \phi^{\beta,c} + \phi^{\alpha,ac} \phi^{\beta,b} \right) + f^{i,\alpha\beta\delta} \phi^{\alpha,a} \phi^{\beta,b} \phi^{\delta,c} \quad [25]$$

where $\alpha, \beta, \in \delta\{1,....,n\}$ and $a,b,c = \{1,...,n\}$ are the indexes for the first, second and third order derivative. $f^{i,\gamma_1...\gamma_p}$ are the partials of the dynamics and are computed as follows:

$$f^{i,\gamma_1...\gamma_p} = \left.\frac{\partial^p f^i(t;\xi_0;t_0)}{\partial \xi_0^{\gamma_1} \cdots \partial \xi_0^{\gamma_p}}\right|_{\xi_0^{\gamma_j}=x_0^{\gamma_j}} \quad [26]$$

Note that the partial derivatives in Eqs. [22] and [26] are calculated with respect to the nominal trajectory. The calculation of the STTs requires the forward integration of $\sum_{q=1}^{s+1} 6^q$ differential equations starting with initial values $\phi_{(t_0,t_0)}^{i,a} = 1$, if $i=a$, and zero otherwise. When the order is $s=3$, the 1554 equations need to be integrated simultaneously. Moreover the computational time and complexity are increased by the numerical evaluations of the analytical partials of the dynamics. In this work, the partials were computed analytically using the symbolic manipulator in the MATLAB(RM) Symbolic Toolbox. As an example, the third order STTs integration along a 2.5 day period required approximately 4 hours using a Windows 7 OS 3.16Ghz Intel[(R)]Core[(TM)]2 Duo CPU.

When implementing the filter it is necessary to calculate the STTs at each intermediate time, in between $t_0$ and t, at which a new measurement is available. The intermediate STTs are called local STTs. Whereas the global STTs map the deviation at the initial time $t_0$ to the deviation at time $t_{k+1}$, the local STTs map the deviation at time $t_k$ to the deviation at time $t_{k+1}$.

There are two methods to compute the local STTs. The first method integrates all the partials from $t_k$ to $t_{k+1}$ without using the information from the global STTs. The second method, indeed, calculates the local STTs from $t_k$ to $t_{k+1}$, having previously integrated the global STTs over the time spans $t_0$-$t_k$ and $t_0$-$t_{k+1}$. The local STTs can be calculated by computing the Inverse State Transition Tensors (ISTTs) via series reversion, as follows:

$$\psi^{i,a} = [\mathbf{\Phi}^{-1}(t,t_0)]^{i,a} \quad [27]$$

$$\psi^{i,ab} = -\psi^{i,\alpha} \phi^{\alpha,j_1,j_2} \psi^{j_1,a} \psi^{j_2,b} \quad [28]$$

$$\psi^{i,abc} = -\left[ \psi^{i,\alpha} \phi^{\alpha,j_1 j_2 j_3} + \psi^{i,\alpha} \left( \phi^{\alpha,j_1} \phi^{\beta,j_2 j_3} + \phi^{\alpha,j_1 j_2} \phi^{\beta,j_3} + \phi^{\alpha,j_1 j_3} \phi^{\beta,j_2} \right) \right] \psi^{j_1,a} \psi^{j_2,b} \psi^{j_3,c} \quad [29]$$

where $j_1, j_2, j_3 = \{1,...,n\}$ are the indexes for the first, second and third order derivative. Note that the series reversion requires the calculation of the inverse of the state transition matrix. The inverse matrix needs to be calculated with high precision otherwise the terms in the expansion result can be affected by a considerable error.

Figure 4 shows the result of the propagation for 2.5 days of the variation vector, with respect to the actual trajectory, $\delta \mathbf{x}_0 = [50\ 50\ 50\ 0.05\ 0.05\ 0.05]^T$ by using third order global STTs, third order integrated local STTs and third order local STTs calculated via series reversion.

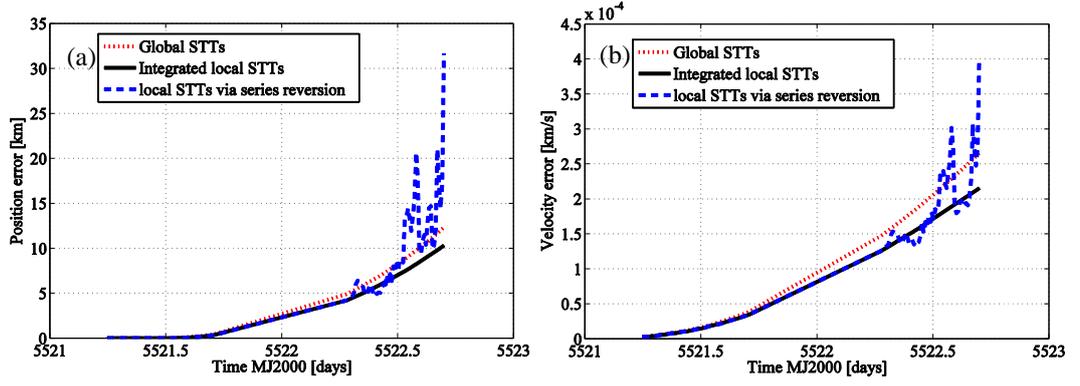

Figure 4: Position (a) and velocity (b) errors using global STTs, integrated local STTs and local STTs obtained via series reversion

The solution calculated using the STTs obtained by series reversion produces oscillations towards the end of the integration period. In order to avoid this problem, in the remainder of this work, the integrated local STTs are used.

Once the state transition tensors are available for some time interval $[t_k, t_{k+1}]$, the mean and covariance matrix of the relative dynamics at $t_k$ can be mapped analytically to $t_{k+1}$ as a function of the probability distribution at $t_k$. From $t_k$ to $t_{k+1}$ the propagated mean and covariance can be stated as:

$$\delta m_{k+1}^i(\delta \mathbf{x}_k) = \sum_{p=1}^{s} \frac{1}{p!} \phi_{(t_{k+1},t_k)}^{i,\gamma_1 \cdots \gamma_p} E[\delta x_k^{\gamma_1} \cdots \delta x_k^{\gamma_p}] \quad [30]$$

$$P_{k+1}^{ij} = E[(\delta x_{k+1}^i - \delta m_{k+1}^i)(\delta x_{k+1}^j - \delta m_{k+1}^j)] = \sum_{p=1}^{s}\sum_{q=1}^{s} \frac{1}{p!q!} \phi_{(t_{k+1},t_k)}^{i,\gamma_1 \cdots \gamma_p} \phi_{(t_{k+1},t_k)}^{j,\varsigma_1 \cdots \varsigma_q} E[\delta x_k^{\gamma_1} \cdots \delta x_k^{\gamma_p} \delta x_k^{\varsigma_1} \cdots \delta x_k^{\varsigma_q}] - \delta m_{k+1}^i \delta m_{k+1}^j \quad [31]$$

where $\{\gamma_i, \varsigma_j\} \in \{1, \ldots, n\}$ are the indexes for the first, second and third order derivative. Figure 5 shows the projection on the x-y plane of a propagated initial distribution with standard deviation $\boldsymbol{\sigma} = [50 \ 50 \ 50 \ 0.05 \ 0.05 \ 0.05]^T$ and zero mean. A total of $10^6$ trajectory samples were propagated forward in time for 2.5 days from TLI using the STTs with $s = \{1, 2, 3\}$ (yellow, green, blue dots) and compared to the Monte Carlo simulation (red dots). The result of the propagation with the first order STTs is an ellipse, while using the second and third order expansions one obtains an approximation more similar to the actual distribution.

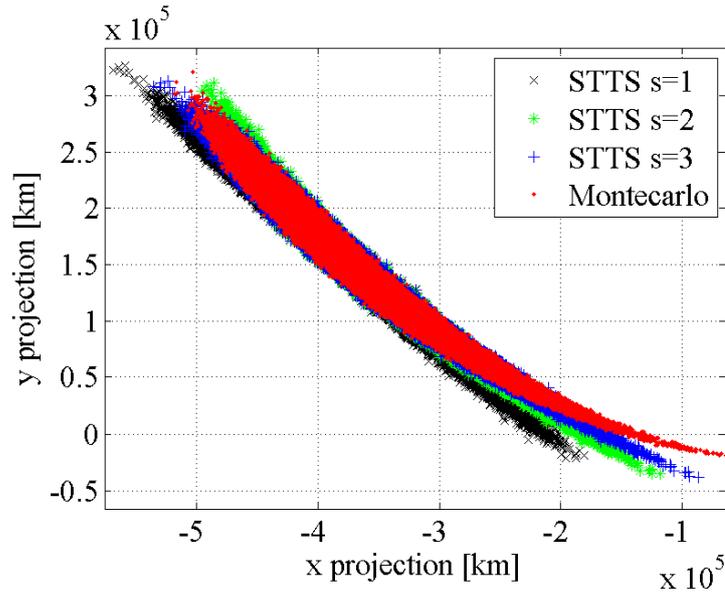

Figure 5: Statistical distribution for STTs and actual Monte Carlo simulations.

Figure 5 suggests that the propagated distribution is not exactly Gaussian though the initial samples were taken from a Gaussian distribution. The assumption in the construction of the filter is, however, that the statistical distribution remains close to Gaussian. This assumption was proved to provide a sufficiently accurate state estimation (see Park and Sheeres 2006b, and Park, 2007). If one sticks to the hypothesis of a Gaussian distribution, the joint characteristic function for a Gaussian random vector can be computed as (Park, 2007)

$$\vartheta(\mathbf{u}) = E[e^{j\mathbf{u}^T\mathbf{x}}] = \exp(j\mathbf{u}^T\mathbf{m} - \frac{1}{2}\mathbf{u}^T\mathbf{P}\mathbf{u}) \qquad [32]$$

where $j = \sqrt{-1}$ and the expected higher moments can be computed using:

$$E[\mathbf{x}^{\gamma_1}\mathbf{x}^{\gamma_2}\cdots\mathbf{x}^{\gamma_p}] = j^{-p} \left.\frac{\partial^p \vartheta(\mathbf{u})}{\partial \mathbf{u}^{\gamma_1}\partial \mathbf{u}^{\gamma_2}\cdots\partial \mathbf{u}^{\gamma_p}}\right|_{\mathbf{u}=0} \qquad [33]$$

then the state prediction and associated covariance can be computed as follows (Park and Sheeres, 2006a):

$$\mathbf{m}^i_{k+1} = \phi^i(t_{k+1};\mathbf{m}^+_k) + \delta\mathbf{m}^i_{k+1} = \phi^i(t_{k+1};\mathbf{m}^+_k) + \sum_{p=1}^{s}\frac{1}{p!}\phi^{i,\gamma_1\cdots\gamma_p}_{(t_{k+1},t_k)} E[\delta x^{\gamma_1}_k \cdots \delta x^{\gamma_p}_k] \qquad [34]$$

$$(\mathbf{P}^-_{k+1})^{ij} = \sum_{p=1}^{s}\sum_{q=1}^{s}\frac{1}{p!q!}\phi^{i,\gamma_1\cdots\gamma_p}_{(t_{k+1},t_k)}\phi^{j,\varsigma_1\cdots\varsigma_q}_{(t_{k+1},t_k)} E[\delta x^{\gamma_1}_k \cdots \delta x^{\gamma_p}_k \delta x^{\varsigma_1}_k \cdots \delta x^{\varsigma_q}_k] - \delta m^i_{k+1}\delta m^j_{k+1} + Q^{ij}_k \qquad [35]$$

For the measurements update phase it is assumed that the linearization of the measurements function provides a sufficient approximation. In this way the Kalman gain $\mathbf{K}_{k+1}$, the state estimate $\mathbf{m}^+_{k+1}$ and the covariance matrix $\mathbf{P}^+_{k+1}$ are analogous to those in Eq. [10] and can be computed as follows:.

$$\mathbf{K}_{k+1} = \mathbf{P}^-_{k+1}\mathbf{H}^T_{k+1}\left[\mathbf{H}_{k+1}\mathbf{P}^-_k\mathbf{H}^T_k + \mathbf{R}\right]^{-1}$$

$$\mathbf{m}^+_{k+1} = \mathbf{m}^-_{k+1} + \mathbf{K}_{k+1}(\mathbf{y}_{k+1} - \mathbf{h}(\mathbf{m}^-_{k+1})) \qquad [35]$$

$$\mathbf{P}^+_{k+1} = (\mathbf{I} - \mathbf{K}_{k+1}\mathbf{H}_{k+1})\mathbf{P}^-_{k+1}$$

where $\mathbf{m}^+_{k+1}$ substitutes $\hat{\mathbf{x}}^+_{k+1}$. The linear assumption simplifies the problem a great deal since the updated phase does not require the computation of the higher order partials of the measurements equations (Park and Sheeres, 2006a). In this way the filter velocity is increased, but at the same time the precision is not affected considerably. Note that when $s=1$, the HAEKF becomes the linear Kalman filter, whose performance is inferior to the EKF, as demonstrated by Mayback, 1982. Since the STTs are integrated offline with respect to the nominal trajectory, the idea is to use second and third order expansions in order to incorporate the nonlinear effects exploiting the advantages given by the pre-integration.

## 4  NAVIGATION

The WSB transfer is conceptually divided into two legs: (1) the first one departs from GTO perigee and extends up to what we call the WSB point. **At** the WSB point nominally a manoeuvre is executed to target the correct return leg that leads to capture by the Moon, (2) the second leg goes from the WSB point to the Moon till the point of Lunar Orbit Insertion.

In order to define the required accuracy for the orbit determination process one can define a capture corridor or region in the state space at every time along the transfer, within which the spacecraft needs to be to achieve the required lunar capture orbit. As it will be shown later in this paper, the capture corridor provides also the basis for the definition of a robust navigation strategy: the idea is to manoeuvre in order to maintain the spacecraft within the capture corridor with sufficient margin to accommodate any orbit determination and navigation error.

The navigation strategy proposed in this paper aims to reduce the deviations with respect to the nominal trajectory, at minimum Δv cost, and guarantee long term capture around the **Moon** in the case of contingencies. Deviations can be due to the inaccuracy in the states at translunar injection, unbalanced accelerations due to the attitude control, the environment (i.e. gravitational perturbation and solar pressure), and errors in the execution of major Δv manoeuvres.

## 4.1 *Lunar Capture Corridor*

In order to define the required accuracy in the estimation of the state of the spacecraft during the transfer to the Moon, the concept of Lunar Capture Corridor is introduced (Zuiani et al., 2011). The capture corridor is the set $B_k$, in the state space, at time $t_{WSB}<t_k<t_{LOI}$ that contains the end state of all the trajectories back propagated towards the WSB region from a ball $B_0$ around the nominal state on the target plane at lunar orbit insertion, without including any correction manoeuvre. The target plane is defined as the plane orthogonal to the velocity for the nominal arrival conditions at LOI.

The states within the ball $B_0$ are such that, after nominal lunar orbit insertion the resulting orbit remains proximal, after a time $T_{LO}$, to the nominal one. Based on the works of Croisard et al., 2009 and Gibbings et al., 2010 the ball $B_0$ was assumed to have maximum deviation of δr=18 km in position and δv= 5 m/s in velocity. A representation of the capture corridor (see Figure 6) was obtained by back-propagating 10000 perturbed state vectors taken within the ball $B_0$ on the target plane at LOI.

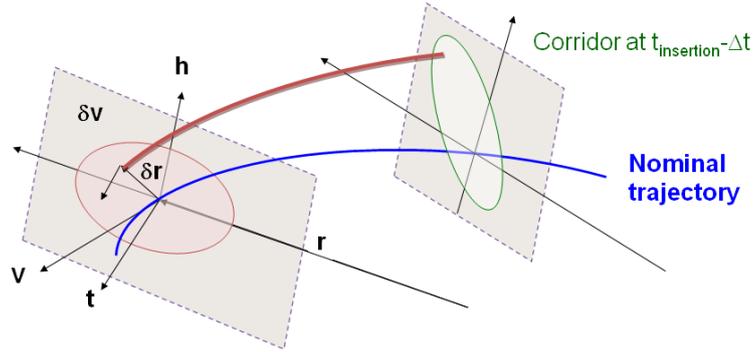

Figure 6: Schematic of the backward propagation in the radial, transversal and out of plane r-t-h reference frame.

Table 2: Orbit parameters of reference lunar operational orbit.

| Parameter | Value |
| --- | --- |
| Semi-major axis | 13084 km |
| Eccentricity | 0.8 |
| Inclination | 56.2° |
| Perigee Altitude | 879 km |
| Apogee Altitude | 21813.2 km |
| Argument of perigee | 270° |
| Ascending node | Free |

The reference lunar orbit (Gibbings et al., 2010) is a highly eccentric frozen orbit with orbital parameters defined in Table 2. The deviated state was established assuming that the post LOI manoeuvre orbit remains proximal to the target reference one. The proximity of the deviated orbit to the nominal one was defined by checking that the deviated lunar orbit remains stable (i.e. does not crash) for a time $T_{LO}$ equal to 6 months, plus or minus 20 days.

**Figure 7 and Figure 8** show the results of the backward propagation (position and velocity dispersion in the *r-h* plane) at two weeks from lunar injection for a typical transfer trajectory. The position and velocity plots report the variation with respect to the nominal value and, therefore, are centered around the origin.

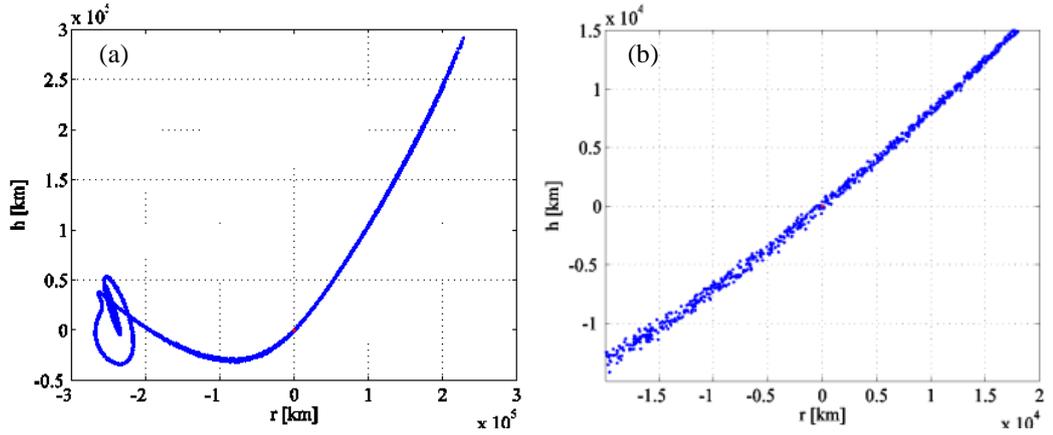

Figure 7: Example of position dispersion in the r-t-h reference system. (a) Large Plot. (b) Close up around the nominal transfer.

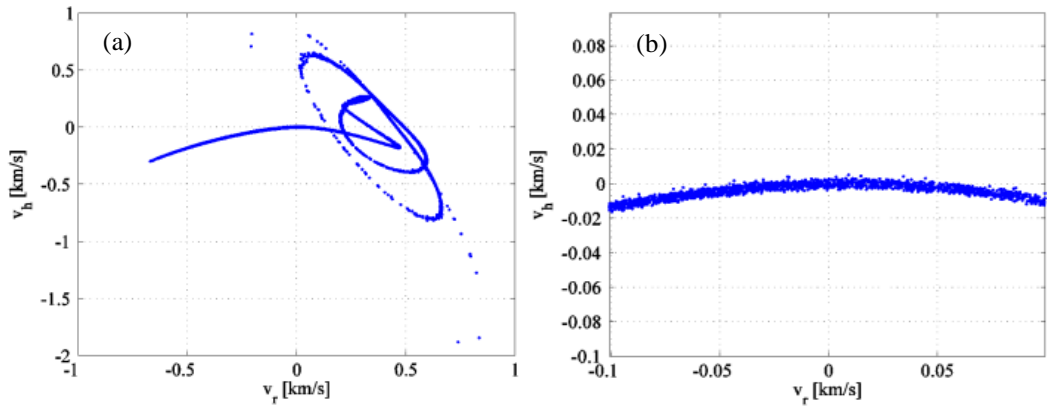

Figure 8: Example of velocity dispersion in the r-t-h reference system. (a) Large plot. (b) Close up around the nominal transfer.

The corridor provides a necessary condition for the accuracy of the orbit determination process. If the orbit determination process is unable to predict that the state of the spacecraft is inside the corridor then it is not possible to know if the satellite will achieve capture. This is of course only a necessary condition, not a sufficient one. The navigation always targets the nominal solution and then it checks that the flow of possible trajectories coming from the errors in state estimation and manoeuvre execution is within the corridor (satisfaction of the necessary condition). The convergence to the required set of states at injection is achieved by retargeting the nominal solution and checking the inclusion in the corridor at different stages along the transfer.

### 4.2 Navigation Strategy

Intermingled along the transfer and in between the orbit determination some Trajectory Correction Manoeuvres (TCMs) are optimally timed and executed to remain within the corridor during the WSB-LOI leg. The goal of each TCM is to minimise the deviation from the nominal trajectory at certain points, called waypoints, along the transfer. The navigation strategy always targets the nominal solution, and then checks that the flow of possible trajectories coming from the errors in state estimation and manoeuvre execution is within the corridor. Therefore, after each orbit determination segment a TCM may, or may not, be required. Orbit determination is assumed to occur over at least a two days period before each planned corrective action. This is to guarantee a good level of convergence of the filter. During the observation period the measurements are assumed to be received every 60 seconds (Thornton and Border, 2003).

Following each TCM, an orbit determination campaign estimates both position and velocity. In this sense, the spacecraft should be visible after and possibly during all TCMs. The sum of all the TCMs will lead to an increase in the mission $\Delta v$ and propellant budget. Two TCMs are allocated after each OD to correct the trajectory up to the next waypoint $t_{\text{wp}}$, i.e. the time along the trajectory where the spacecraft is expected to be within the corridor. The waypoints are equally spaced in time. The last waypoint coincides with the nominal arrival conditions at LOI. After the first correction manoeuvre, another orbit determination takes place and the operation is repeated in order to redefine the following manoeuvre and check if the first correction has been performed correctly. At each waypoint, the nominal state of the spacecraft is $\mathbf{x}_{nominal}(t_{wp})$ and the state provided by the implementation of the TCMs is $\mathbf{x}(t_{wp})$. Each TCM is defined by its time of execution $t_{TCM}$ and the components of the velocity variation with respect to the local velocity vector. The following constrained optimisation problem is then solved to optimally allocate and size each pair of TCMs:

$$\min_{\mathbf{u} \in U} \Delta v_{TCM_1} + \Delta v_{TCM_2}$$
$$s.t. \quad [36]$$
$$\mathbf{x}(t_{wp}, \mathbf{u}) - \mathbf{x}_{nominal}(t_{wp}) = 0$$

where the control vector $\mathbf{u}$ contains the TCM's time of execution and $\Delta v$ components. Thus, the scheduling (time and date), direction and magnitude of the **TCM** must be optimised. Problem [36] is solved by using the function *fmincon* of the 2010 MATLAB$^{RM}$ Optimisation Toolbox).

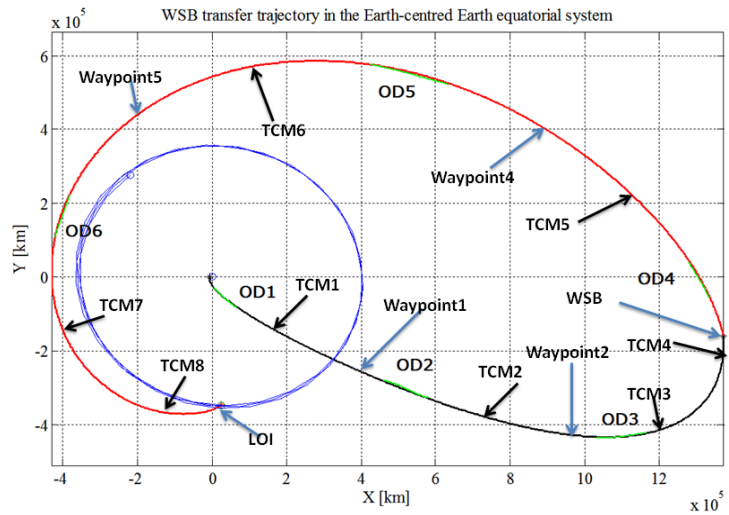

Figure 9: Navigation Strategy along the WSB transfer.

As shown in Figure 9, during the first and second leg the last waypoint coincides with the nominal orbit state vector at WSB (after the WSB correction) and LOI (before the manoeuvre) respectively.

### 4.2.1 Orbit determination effects

For the case of ESMO, a single ground station is used for orbit determination. However, the ground station, in Raisting, Germany, is currently still under development. Thus, it is not possible to rely on a defined level of accuracy when dealing with the measurements. Especially angular measurements are uncertain, although the ground segment is currently planning to achieve a 10 mdeg accuracy. A sensitivity analysis is therefore required to establish which level of accuracy is acceptable to navigate to the Moon. The measurement error (1σ) in range and range-rate was assumed to be between 10 and 20 m and between 0.5 and 1 mm/s respectively. These are conservative values according to the current state of the art (Thornton and Border, 2003) for orbit determination. The error in the angle estimation was assumed to be between 3 and 80 mdeg. Then, 100 full trajectories were simulated to assess the effect of errors in range, range-rate and angles. Table 3 reports the errors that were used in this analysis.

Table 3 : Estimated measurement errors.

|  | Error | | |
|---|---|---|---|
| Range [m] (*1σ*) | 10 | 15 | 20 |
| Range-rate [mm/s] (*1σ*) | 0.5 | 0.75 | 1 |
| Angles [mdeg] (*1σ*) | 3 | 20 | 80 |

Note that the lower limit for range and range-rate are consistent with Landgraf, 2007, while the upper limit for range has been set equal to the 1992 accuracy in Thorton and Border, 2003. Finally the angles upper limit is set equal to the accuracy given by Villafranca GS (Spain). Figure 10 reports the sensitivity analysis for a 5465 MJD2000 launch. Manoeuvres are affected by 1% errors. Each dot corresponds to the orbit determination accuracy for a given combination of the nine errors in Table 3. Colours are associated to the magnitude of the velocity error at LOI (Figure 10a), position error at LOI (Figure 10b) and navigation budget (Figure 10c). The value in the figures is the mean plus the variance over all the runs. Red colour is associated with the maximum value of the scale, while blue is associated with the smallest value of the colour scale.

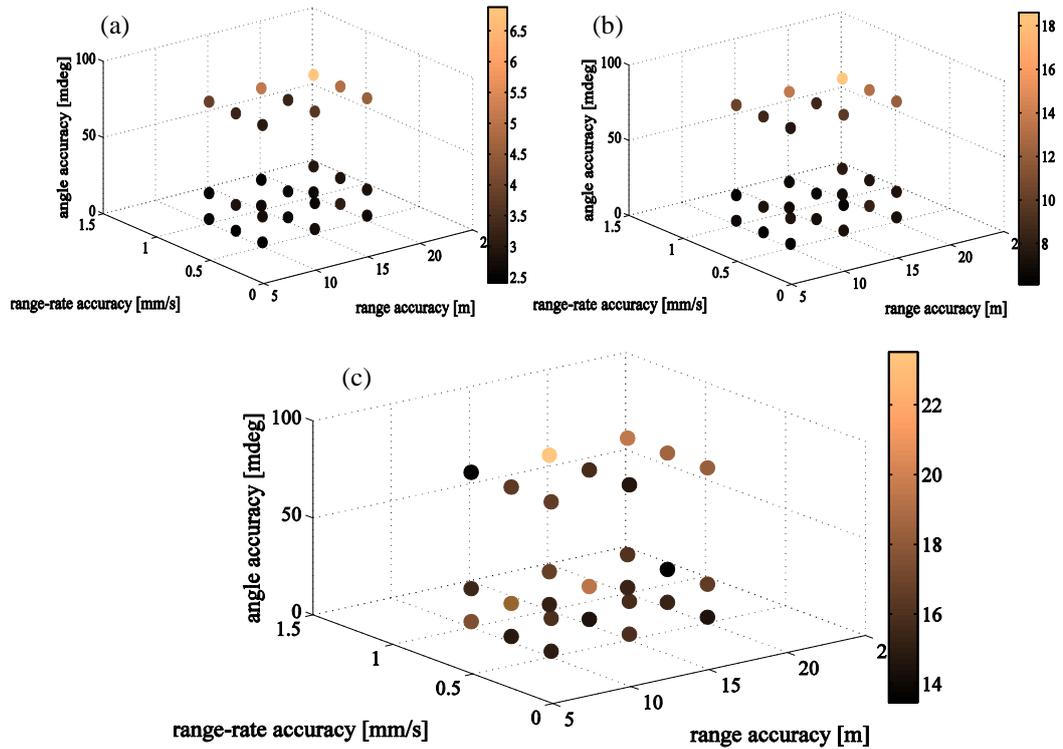

Figure 10: Sensitivity analysis for range, range-rate and angles for 1% TCM error: (a) velocity error at LOI, (b) position error at LOI and (c) navigation budget.

This analysis shows that in general the performance is good when the angles accuracy is up to 20 mdeg. The accuracy at LOI allows the spacecraft to be within the boundaries defined by the corridor while keeping the navigation budget lower than 25 m/s. The 20 mdeg accuracy guarantees good insertion accuracy, and reduces the effects on the navigation budget. In the following sections analyses will consider OD measurement errors (1σ) of 15 m on range, of 0.5 mm/s on range-rate and 20 mdeg on angles.

### 4.2.2 Filter Comparison

After defining the measurements accuracy in the previous section, this section presents a comparison among the three filters described in Section 3.1 in order to identify the most suitable filtering technique. The impact of the different methods on the navigation will be addressed in Section 4.3.2.

As a representative example, this section reports the results for the first orbit determination campaign, which is allocated soon after the final translunar injection manoeuvre. It is considered that a 2.5 days orbit determination campaign is performed after the translunar injection manoeuvre. The set of pseudo-measurements is taken every 10 minutes, when the spacecraft is visible from the ground station. The time step between measurement samples is higher than the one used in the following sections to emphasize the effects of dynamics nonlinearities. In fact, a desirable benefit of high order nonlinear filters would be to reduce the number of steps to obtain the same predication accuracy. It is therefore important that the reduction of the number of steps compensates for the higher cost of each step.

During the filtering process comparison, then, the same measurements have been used and it is assumed that the system dynamics is known completely.

The dispersion of the initial states of the spacecraft is assumed to follow a Gaussian distribution with zero mean and standard deviation equal to 100 km in position and 0.1 km/s in velocity similarly to Park and Sheeres, 2006b. The initial covariance matrix is diagonal with each component along the diagonal equal to the squared of the standard deviation of position and velocity. The mean and covariance matrix are mapped between each step $k$ using the STTs for $s = \{1, 2, 3\}$, see Eq. [21], the first order is used to for comparison against the linear propagation in the EKF and the second and third order are used to compare the HAEKF against the unscented transformation in the UKF. A Monte Carlo simulation based on $10^6$ samples is then used to validate the correctness of the outcome of each filter. In the third case the calculation of the mean and covariance matrix is as follows:

$$\mathbf{m}_i(t) = \frac{1}{N} \sum_{k=1}^{N} \phi_i(t; x_k^0, t^0) \quad [37]$$

$$\mathbf{P}_{ij}(t) = \frac{1}{N-1} \sum_{k=1}^{N} [\mathbf{m}_i(t) - \phi_i(t; x_k^0, t^0)][\sum_{k=1}^{N}[\mathbf{m}_j(t) - \phi_j(t; x_k^0, t^0)] \quad [38]$$

Figure 11 shows the mean and the projection of 1-σ covariance matrix onto the x-y plane after 2.5 days. Red dots represent all the outcomes form the Monte Carlo simulation.

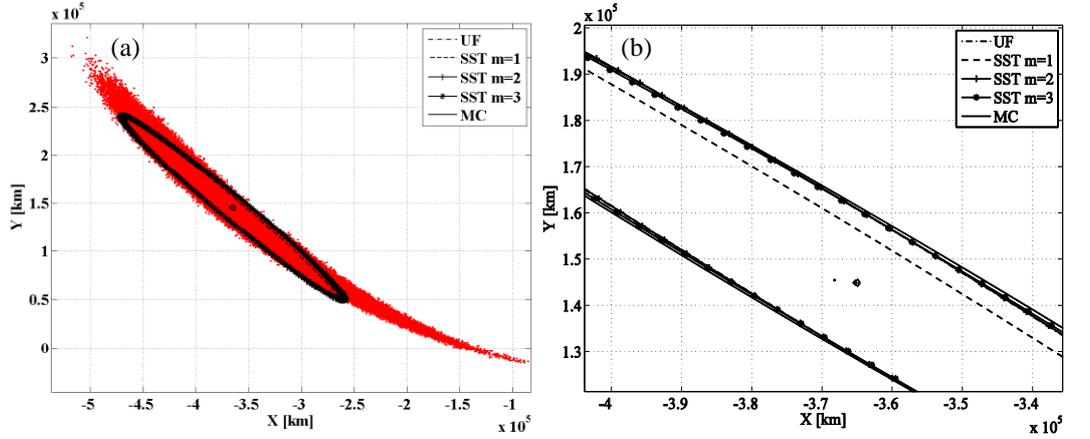

Figure 11: Mean and 1-σ error ellipsoid projected onto the position x-y plane after being propagated for 2.5 days: (a) large plot, (b) close up around the mean value.

Assuming that the Monte Carlo simulation is the true solution, the third order solution is the most accurate approximation, and the second order and unscented transformation are approximately the same, thus confirming that UT is able to propagate the covariance matrix with the precision equivalent to a second order expansion. Table 4 reports the mean position and velocity error and the maximum square root diagonal element difference with respect to the Monte Carlo solution. Expectedly, the third error propagation is the most accurate for the covariance matrix, while first order solution is the worst, with a significant error also in the mean value.

Table 4: Mean and covariance matrix simulation using the different techniques. CPU time using a 64-bit Linux CentOS 5.4 2.93GHz Intel(R) Xeon(TM) X5570.

|  | UT | STTs s=1 | STTs s=2 | STTs s=3 |
|---|---|---|---|---|
| CPU time [s] | 1.234 | 1.105 | 1009.786 | 12429.132 |
| Mean position error [km] | 194.72 | 3370.32 | 195.30 | 195.30 |
| Mean velocity error [km/s] | 0.00332 | 0.03477 | 0.00333 | 0.0333 |
| Maximum square root diagonal element difference | 10691 | 11591 | 10575 | 4348 |

The second, third order expansion and the unscented transformation are able to propagate the mean orbit with a relative error of approximately 200 km. Note that in this case the mean initial deviation is zero, therefore the third order for the expected value is equal to zero, and the third order mean coincides with the one from the second order expansion. The required CPU times vary considerably from the first to the third order approximation. The CPU time includes all the operations necessary to integrate STTs through Eqs. [23]–[25] and calculate propagated mean and covariance matrix using Eqs. [34]-[35] (without including the contribution of the matrix $Q$ there reported). The propagation with third order STTs required about 27% the time required to launch $10^6$ Monte Carlo simulations, which took about 45275 s.

During the filtering process comparison, the initial guess was set equal to the state vector of the nominal unperturbed trajectory at the time the first measurement is received. In the evaluation of the three filters, many different sets of pseudo measurements were simulated but the results showed negligible difference in the filter performance given the state and measurement uncertainties considered in the comparison.

It has to be underlined that the HAEKF approach required the integration of 258, for $s=2$, and 1554, for $s=3$, differential equations with respect to the reference unperturbed trajectory. This integration was performed prior to using the STTs in the actual filter process but it was done for each time span corresponding to the measurement update. The dimensional dynamic equations were numerically integrated with an explicit, variable step size, Runge-Kutta integration method with a $10^{-9}$ and $10^{-9}$ relative and absolute accuracies respectively.

As Figure 12 shows, the second and third order HAEKF and UKF present better convergence results and produce a more accurate estimate than the EKF. The increments in velocity and position error (straight line in Figure 12) are due to the propagation of the last estimate since the spacecraft is not visible during that period and measurements are not available.

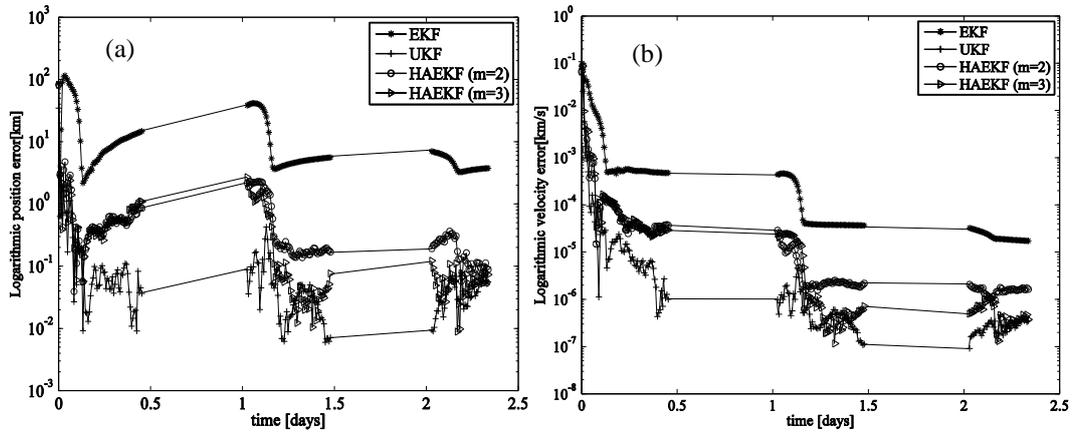

Figure 12: Comparison of the absolute errors for 100 km position uncertainties (a) and 0.1 km/s velocity uncertainties (b).

If an initial uncertainty in position and velocity equal to 1% of the nominal state vector is considered, the UKF shows superior performance over the other filter methods, as shown in Figure 13. The absolute errors are computed more accurately by the UKF, with a final position error lower than 1 km and velocity error lower than $10^{-5}$ km/s. The higher precision is due to the fact that the measurements model is used without introducing any linearization assumption.

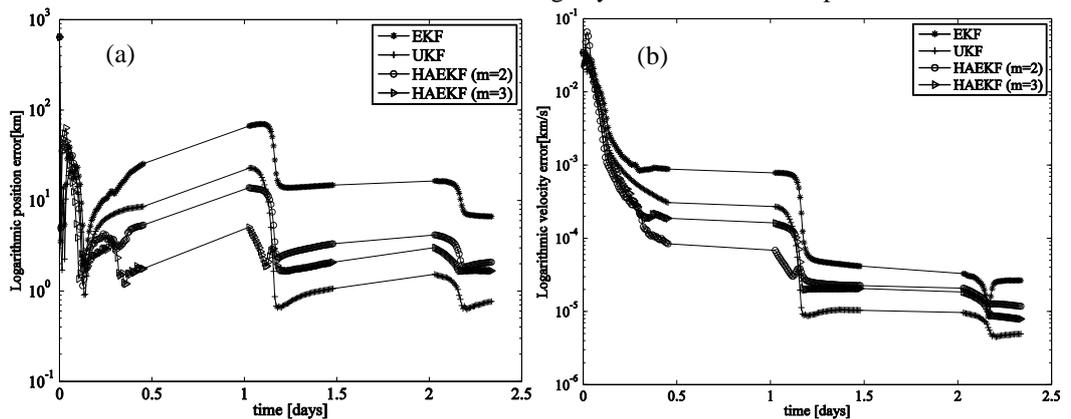

Figure 13: Comparison of the absolute errors for uncertainties on position (a) and velocity (b) equal to 1% of the nominal unperturbed trajectory.

Even if the third order HAEKF provides a solution essentially equivalent to the UKF when the initial guess is relatively close to the actual position and good results in the other case, it yields a much slower filter process. Table 5 reports the elapsed CPU time necessary to process all the measurements by the different filters during the 2.5 days observational period.

Table 5: Elapsed CPU time for the filtering processes, using a 64-bit Linux CentOS 5.4 2.93GHz Intel[(R)] Xeon[(TM)] X5570, and absolute error for the estimated state.

|  | EKF | UKF | HAEKF s=2 | HAEKF s=3 |
|---|---|---|---|---|
| Elapsed CPU time [s] | 20.52 | 47.02 | 45.06 | 1737.22 |
| Absolute position error [km] | 6.657 | 0.765 | 2.081 | 1.602 |
| Absolute velocity error [km/s] | $2.648 \cdot 10^{-5}$ | $4.939 \cdot 10^{-6}$ | $1.117 \cdot 10^{-5}$ | $7.905 \cdot 10^{-6}$ |

The most computational intensive filter is represented by the third order HAEKF itself, while the elapsed CPU time for the second order HAEKF and UKF is about 50 seconds. The EKF is the fastest method but it is still the least precise when compared to the other methods.

The EKF required the integration of $n + n^2 = 42$ equations (necessary for the update of the estimated state and the integration of the STM), being $n = 6$, and the UKF required the integration of $(2n+1)n = 78$ equations (for the updated estimated state and the propagation of the sigma points) between each measurement update. Moreover, the higher cost of the UKF with respect to the EKF is also given by the definition of the innovation process in Eq. [18] and Eq. [19]. Although the HAEKF required only 6 integrations between each measurement update in Eq.[34], less than the UKF and EKF, the number of joint functions evaluations and summations required to compute the mean value and covariance matrix in Eq. [34]-[35] between each measurement update is $\sum_{q=1}^{s} 6^q + \left( \sum_{p=1}^{s+1} 6^{p+1} \right)^2$, which is equal to 66822 and 2397854 for $s = 2$ and $s = 3$ respectively.

### 4.2.3 Impact on Navigation

The major impact of a correct estimation of the state of the spacecraft is on navigation. A poor estimate could lead to an incorrect implementation of TCM. In order to compare the impact of the three filters on the navigation strategy, problem [36] was solved using the state estimate coming from the EKF and UK, while for the HAEKF the following problem proposed by Park and Sheeres, 2006b was solved:

$$\min_{\Delta \mathbf{v}_1^{TCM}} \left( \delta \mathbf{m}_r^{waypoint} \right)^2 = \left( E[\phi_r(t^{waypoint}; \delta \mathbf{m}_r^{TCM}, \delta \mathbf{m}_v^{TCM} + \Delta \mathbf{v}_1^{TCM}; t^{TCM})] \right)^2 \quad [39]$$

with:

$$\Delta \mathbf{v}_2^{TCM} = -\delta \mathbf{m}_v^{waypoint} \quad [40]$$

This method represents the so called nonlinear statistical targeting correction method introduced by Park and Sheeres, 2006b for the high-order expansion semi-analytical method. By using Eq.[30], the statistical information on the variance of the state estimated at the end of the OD, is included in the propagation of the mean trajectory. Thus the correction manoeuvre is based on the average of all the possible the propagated trajectories and not on the propagation of the estimated mean state. In this way the TCMs are setting to zero the mean of all trajectory deviations at the waypoint. Eq. [39] is solved using *fminsearch* of the 2010 MATLAB[RM] Optimisation Toolbox.

It has to be pointed out that if the manoeuvres are performed, the actual state error at the waypoint will be not zero, but the statistical mean error will be zero, if the STTs order is sufficiently accurate to take into account all the nonlinear effects. In this comparison, the time and location of the TCM and waypoint was fixed for all the three filters and to be consistent with problem [39], problem [36] was modified by imposing a condition on the final velocity equivalent to [40].

Table 6 reports the results for both the optimized manoeuvre and statistical targeting correction applied to the first waypoint. Since the orbit determination estimate for the UKF proved to be the most accurate, optimised manoeuvres were calculated on the basis of 100 runs, while the nonlinear-statistical targeting corrections were calculated for the HAEKF ($s = 2, 3$) results. The actual trajectory is the one corresponding to the propagation of 1% perturbed trajectory whose estimated state is reported in Figure 13. The waypoint was placed at six days after the end of the orbit determination. It is assumed that the correction manoeuvres are not affected by errors. As it can be seen, the results are almost the same in terms of total correction for all but EKF filters. This is achieved due to the fact the effects of the nonlinearity are included in both methods. The best actual state is achieved by optimising the TCMs on the basis of the estimated state by UKF. The optimized manoeuvre method allows us to consider the actual nonlinear dynamics, while higher orders introduce an approximation. The CPU time difference is considerable for the optimised manoeuvre. In this case the STTs approach results faster because the trajectory is not integrated, and the computational cost is due to the *fminsearch* operations.

Table 6: OD impacts on Navigation using optimized and nonlinear statistical targeting correction methods. A 64-bit Linux CentOS 5.4 2.93GHz Intel$^{(R)}$ Xeon$^{(TM)}$ X5570 was used.

|  | Optimised manoeuvre correction | Statistical targeting correction | |
|---|---|---|---|
|  | UKF | HAEKF s=2 | HAEKF s=3 |
| CPU time [s] | 2452.721 | 0.599 | 27.974 |
| Total correction Δv [m/s] | 77.104 | 77.083 | 77.101 |
| Actual position error [km] | 2.89 | 12.10 | 10.56 |
| Actual velocity error [m/s] | 0.003 | 0.032 | 0.025 |

As it can be seen in Table 6, the results are almost the same in terms of total correction for all but EKF filters.

The optimized manoeuvre method allows us to consider the actual nonlinear dynamics, while higher orders introduce an approximation. The CPU time difference is considerable for the optimised manoeuvre. In this case the STTs approach results faster because the trajectory is not integrated, and the computational cost is due to the *fminsearch* operations.

As pointed out by Park and Sheeres, 2006b, when the navigation data are accurate, both the correction manoeuvres are essentially the same. The nonlinear statistical targeting correction depends on the statistical knowledge about the state vector (given by the covariance matrix), whereas the optimised manoeuvre relies only on the mean value.

The statistical targeting correction method results to be a flexible method since the evaluation of the manoeuvre can be done in one minimisation, while the optimised method need to be assessed by a Monte Carlo simulation. Anyway in the case of the statistical targeting correction method the effects of nonlinearities are significant and it would be necessary to increase the STTs order since the final actual state is less precise than the one obtained using the optimised method.

The unscented Kalman filter was adopted as the baseline filter for the mission since it proved to be fast and able to provide better estimates than the other studied method. Moreover it allowed implementing the adopted navigation strategy, reducing the effect of the navigation errors.

## 5   NAVIGATION RESULTS

This section reports that results of a navigation analysis for transfers to the Moon in the launch window 2014-2015. A Monte Carlo analysis was performed by running 100 simulations of the combined OD and navigation process along each trajectory from departure from GTO to lunar orbit insertion. For this kind of estimation 100 samples corresponds to a confidence interval of 6% with a confidence level of 92%. All typical sources of error both on measurements and manoeuvre execution were included in the simulations. Figures 14 and 15 show the error in position and velocity at LOI and the navigation budget for a representative sample of transfer trajectories. In particular, Figure 14 reports the error at LOI and the navigation budget for a 1% error in the execution of the TCMs. Figures 15 reports the same results but for a 3% error in the execution of the TCMs. The 1-σ error in the measurements for all simulations is: range 15 m, range-rate 0.5 mm/s, angles 20 mdeg.

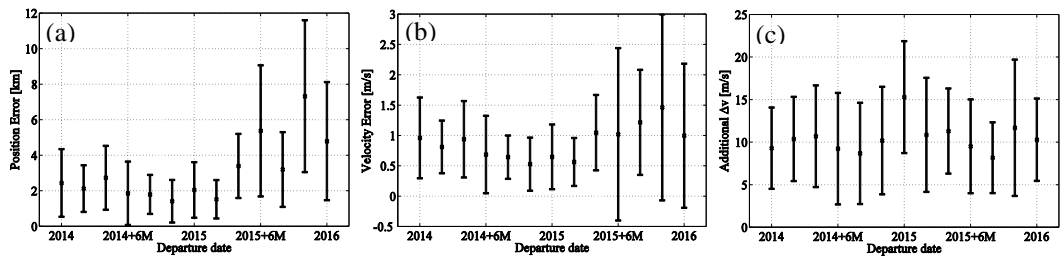

Figure 14: Navigation statistics for 2014-2015 launch window for 1% manoeuvre error: (a) position error at LOI, (b) velocity error at LOI, (c) navigation budget.

The dots inside the bars represent the mean values while the distance between the dot and the upper or lower end is the standard deviation. Figures 14(a) and 14(b) show that ESMO is able to be inserted correctly at LOI according to the adopted navigation strategy, with a low navigation budget (Figure 14(c)). On the contrary when the error affecting TCMs is higher, the accuracy of the orbit insertion is reduced. Figures 15(a) and 15(b) show that in some cases the spacecraft cannot be guaranteed to be within the lunar corridor at LOI. Figure 15(c) show that the navigation budget is increased with many trajectories requiring a higher than 50 m/s in the case of 3% TCM error. Analyses were also conducted for 5% TCM errors, for gives a significantly larger error in position and velocity and a navigation budget that increases up to 100 m/s.

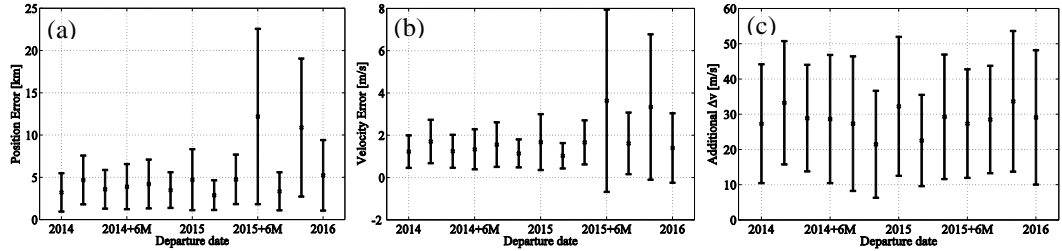

Figure 15: Navigation statistics for 2014-2015 launch window for 3% manoeuvre error: (a) position error at LOI, (b) velocity error at LOI, (c) navigation budget.

A first conclusion, therefore, is that the control of the thrust vector during TCMs deeply affects the navigation budget, and the possibility to achieve the correct lunar insertion LOI point. This is understandable as this is a closed loop control problem and the error in the TCMs reduces the controllability of the system.

## 6   WEAK CAPTURE ANALYSIS

The possibility to achieve a stable lunar orbit, even in the presence of a total failure of the main engine is addressed in this section. Such an event is assumed to occur some days before the nominal LOI. The goal is to achieve a sufficient lifetime of the orbit around the Moon even without the lunar orbit insertion manoeuvre. The lifetime of the spacecraft around the Moon is defined as the difference between the instant at which the nominal LOI manoeuvre should occur and the instant at which the spacecraft impacts the lunar surface or leaves a closed lunar orbit whose maximum distance from the Moon is as high as 60,000 km. In general, as explained in Section 2.1.1. one can assume that when the 2-body energy of the spacecraft with respect to the Moon is non-positive at LOI, the spacecraft can theoretically be captured by the Moon into a weakly stable lunar orbit for a limited period of time. Figure 16 reports the nominal lifetime with respect to the $C_3$ energy before LOI manoeuvre for a 18000 trajectories sample. Many negative $C_3$ trajectories present an uncontrolled lifetime longer than 20 days.

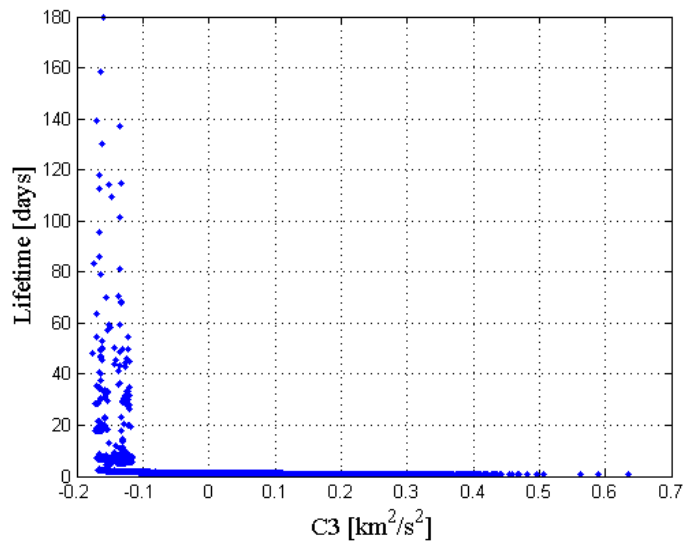

Figure 16: Lifetime vs. $C_3$ for a sample of 18000 nominal transfer trajectories propagated for 180 days after the arrival in LOI. No LOI manoeuvre has been performed.

This allows the attitude engines to be used in order to manoeuvre the spacecraft and extend the spacecraft's lifetime. Figure 17 shows how the trajectory changes considerably if the LOI manoeuvre is not performed. In this case the spacecraft impacts the Moon surface after 32 days since nominal LOI condition.

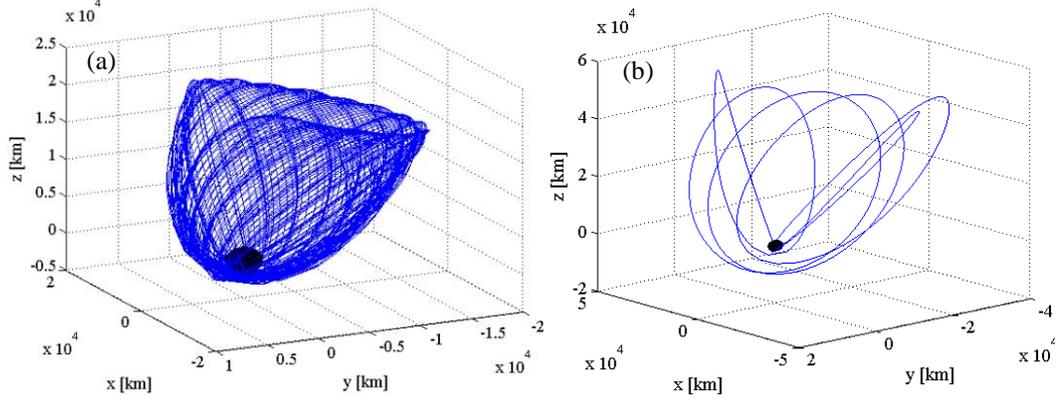

Figure 17: Spacecraft trajectory around the Moon: nominal trajectory after LOI manoeuvre (a) and trajectory in the case of contingency (b).

Since the propellant available for attitude control at lunar arrival will be limited, the analyses were performed considering only 30% of the total attitude control $\Delta v$ to be available for orbit correction post nominal LOI. The total attitude budget for the ESMO mission is about 10 m/s, thus only 3 m/s are available for orbit correction. Since the goal is to extend the lifetime of the orbit, the following constrained optimisation problem was solved with *fmincon*:

$$\min_{\mathbf{u}\in U}(t_{lifetime}^{nominal} - t_{lifetime}^{actual}(\mathbf{u}))$$
$$s.t. \qquad [40]$$
$$\Delta v - \Delta v_{available} = 0$$

The minimization problem seeks for a correction, which is optimally allocated after the nominal arrival time at LOI. The manoeuvre uses the available propellant, in order to minimize the difference between the nominal lunar orbit lifetime $t_{lifetime}^{nominal}$ (required to be 180 days) and the actual lunar orbit after the correction $t_{lifetime}^{actual}$, which means to maximize the lifetime up to 180 days. A correction manoeuvre is optimally allocated after the nominal arrival time at LOI.

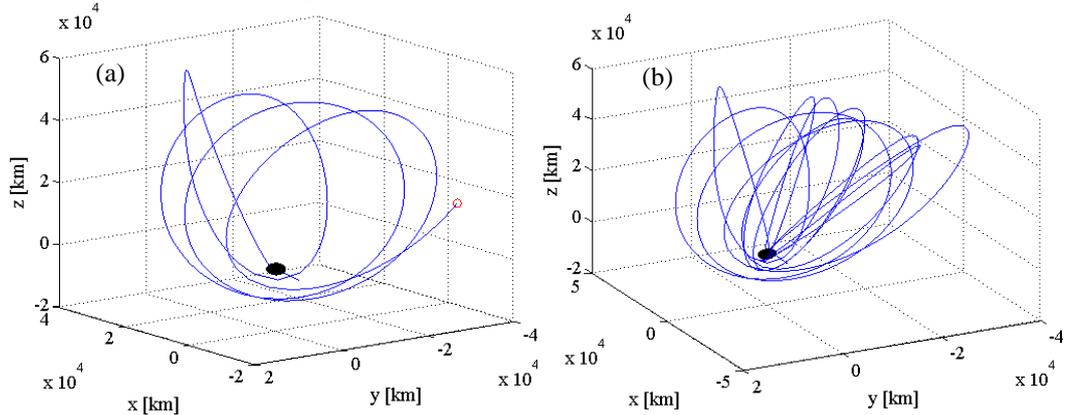

Figure 18: Weak capture: Correction manoeuvre is performed (round point (a)) and corrected trajectory (b).

Figure 18 shows that after a 3 m/s correction manoeuvre is performed (left figure), the number of revolutions is increased thus extending lifetime.

### 6.1 *Weak Capture and Capture Corridor*

The navigation strategy during the WSB-LOI leg aims at maintaining the spacecraft within the capture corridor with sufficient margin to accommodate any orbit determination error. For this reason, the inclusion of the trajectory inside the corridor at the end of the orbit determination is checked together with the estimated state at LOI. If one of these conditions is not met, TCMs are allocated. Figure 19 reports the position and velocity dispersion in the r-t plane for 100 navigation simulations of the reference trajectory, after the last TCM. The TCM is allocated after the last orbit determination at about 4 days before LOI. If the TCMs accuracy is in the order of 1%, all the 100 trajectories are within the corridor soon after the last trajectory correction manoeuvre.

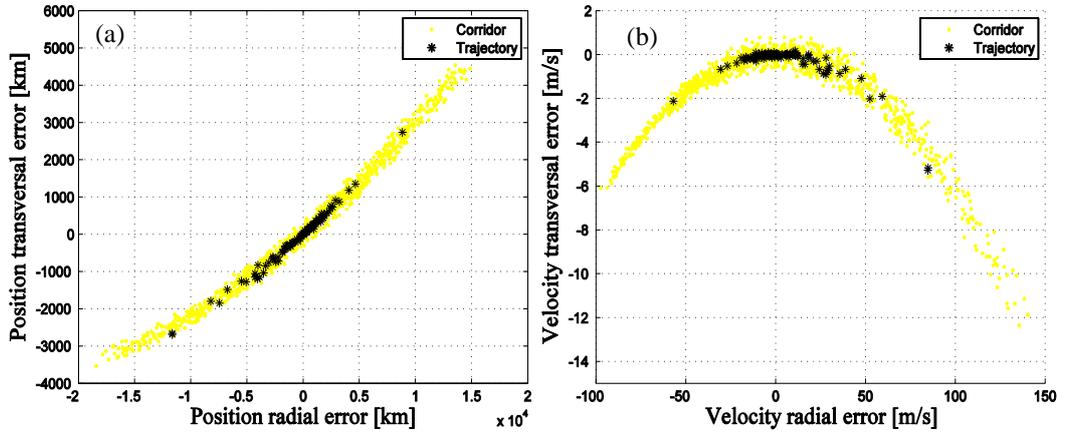

Figure 19: Position (a) and velocity (b) dispersion error for 100 runs after the last TCM (4 days before LOI) in the r-t plane. The yellow points represent the corridor. 1% error on TCM. Nominal arrival time at LOI: 5501 MJD2000. $C_3$=-0.13355 km$^2$/s$^2$.

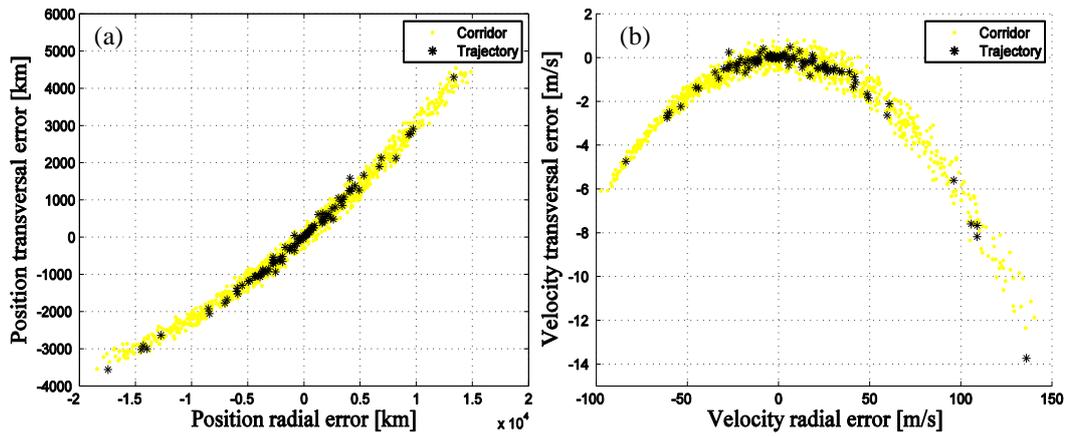

Figure 20: Position and velocity dispersion error for 100 runs after the last TCM (4 days before LOI) in the r-t plane. The yellow points represent the corridor. 3% error on TCM. Nominal arrival time at LOI: 5501 MJD2000. $C_3$=-0.13355 km$^2$/s$^2$.

Figure 20 shows that many trajectories result to be out of the corridor after the last TCM when the error in the execution of the TCMs is 3% or higher. If now one assumes that TCM are affected by a 1% error, the analysis of the capture corridor can help to understand the effect of the solution of problem [39].

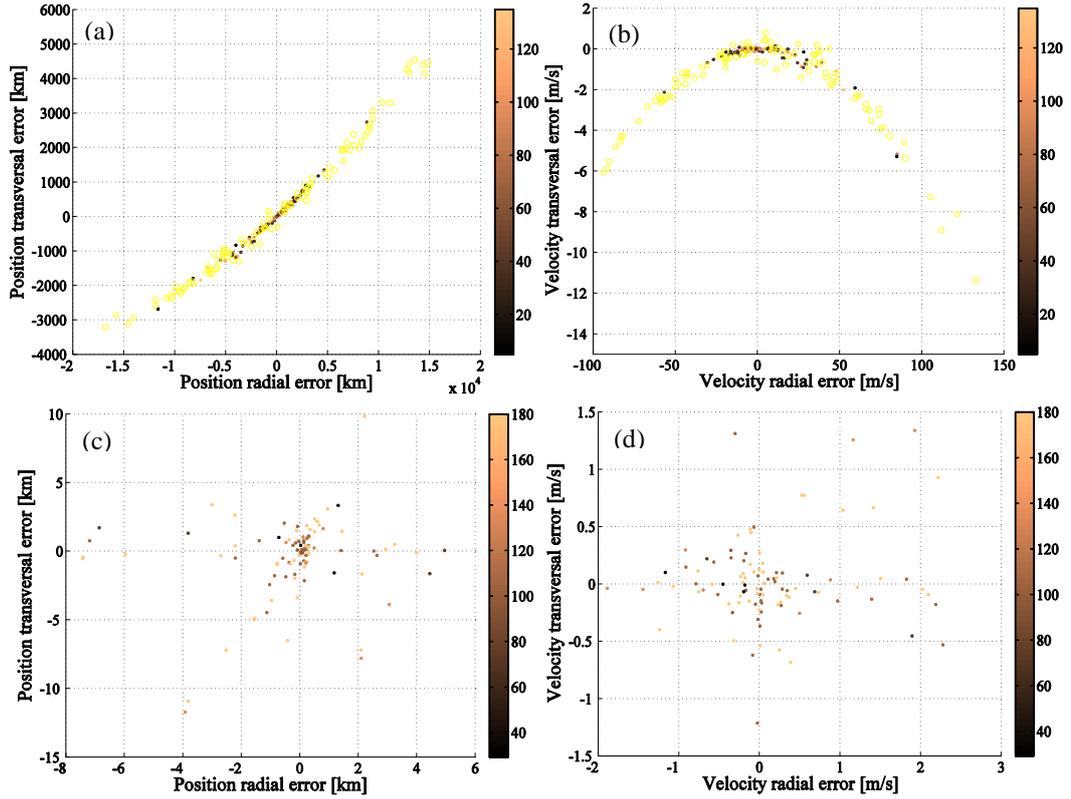

Figure 21: Spacecraft lifetime after around the Moon with respect to position and velocity dispersion in the r-t plane before (a),(b) and after (c),(d) the correction manoeuvre. 1% error on TCMs. Nominal arrival time at LOI: 5501 MJD2000. $C_3$=-0.13355 km$^2$/s$^2$.

As shown in Figure 16, some negative $C_3$ energy trajectories present the possibility to be weakly captured by the Moon for more than 20 days, even if the main engine fails and the LOI manoeuvre is not performed. This condition depends not only on $C_3$ level but also on the relative position of the Sun-Earth-Moon when the spacecraft reaches the LOI. If problem [40] is solved a small correction manoeuvre can be allocated to significantly extend the lifetime of the nominal orbit. Figure 21 show the orbit lifetime for each black point in Figure 19 assuming that a 3 m/s correction manoeuvre is performed after the nominal LOI. The colour is proportional to the lifetime: red is associated to the longest lifetime and blue is the shortest. For the 100 simulations of Figure 19 the mean lifetime increases from about 71 to 138 days. Note that a small error difference can correspond to a high lifetime difference, due to both perturbations and nonlinearities which characterize the system dynamics. Table 7 reports the mean orbit lifetime for 4 negative $C_3$ trajectories considering 100 simulations for each reference trajectory.

Table 7: Contingency analysis for 4 different trajectories. Mean lifetime before and after the correction manoeuvre.

| Trajectory | | | Lifetime before and after TCM(days) | | |
|---|---|---|---|---|---|
| MJ2000 (days) | C3 (km$^2$/s$^2$) | | 1% | 3% | 5% |
| 5501 | -0.13355 | before | 70.8 | 71.1 | 63.8 |
| | | after | 138.4 | 75.5 | 69.4 |
| 5636 | -0.13215 | before | 50.5 | 40.4 | 30.1 |
| | | after | 111.6 | 89.1 | 45.4 |
| 5639 | -0.16323 | before | 60.1 | 58.2 | 59.1 |
| | | after | 99.2 | 91.9 | 67.4 |
| 5606 | -0.16938 | before | 92.1 | 107.1 | 106 |
| | | after | 174 | 124.6 | 171 |

The examples in the table includes also the case in which the error in the execution of the TCMs is up to 5%. The first three rows in the table are cases in which some trajectories are outside the corridor if a 3% and 5% TCM error is considered, while for the last row all the simulated trajectories were within the corridor.

## 7 CONCLUDING REMARKS

This paper presented a strategy to navigate to the Moon along low-energy transfers. The European Student Moon Orbiter was taken as case study. Typical uncertainties in orbit determination, dynamics and manoeuvre execution were included in the analysis. Three types of Kalman filter were studied in order to assess the one able to perform a fast and precise trajectory estimate. The EKF, UKF and HAEKF filters were compared using the same initial conditions. The result of the comparison showed that all the filters are suitable, but best results are obtaining for third order HAEKF and for the UKF. However, the computational cost of the filters led to choose the UKF. On the other hand the stochastic targeting appears to be advantageous from a computational point of view if the propagation of the partials can be performed a priori. It is therefore possible to envisage a combination of state estimation using the UKF and a correction planning using high order expansions.

The contingency analysis at the Moon has shown that it is possible to extend the lifetime of the lunar orbit with a small manoeuvre performed with the attitude engines. An interesting result is that by keeping the spacecraft inside the corridor at the lunar orbit insertion, the possibility to manoeuvre using attitude engines and achieve a stable uncontrolled orbit around the Moon is increased, even if the LOI manoeuvre cannot be performed. Further analysis will also consider the use of different level of available propellant.